\documentclass [12pt,a4paper]{amsart}% Language setting
\usepackage[english]{babel}
\usepackage[letterpaper,top=2cm,bottom=2cm,left=3cm,right=3cm,marginparwidth=1.75cm]{geometry}
\usepackage{amsmath, amsfonts, amssymb, amscd, bezier, amsthm, mathrsfs}
\usepackage{graphicx}
\usepackage{hyperref}
\usepackage[utf8]{inputenc}
\usepackage[T1]{fontenc}
\usepackage{float}
\usepackage[euler]{textgreek}
\usepackage{url}
\usepackage{subcaption}

\newtheorem{theorem}{Theorem}[section]
\newtheorem{lemma}[theorem]{Lemma}
\newtheorem{proposition}[theorem]{Proposition}
\newtheorem{corollary}[theorem]{Corollary}
\newtheorem{question}{Question}

\theoremstyle{definition}

\theoremstyle{remark}
\newtheorem*{remark}{Remark}
\usepackage[pdftex]{color}

\newcommand{\Hyp}{\mathbb{H}}

\newcommand{\Q}{\mathbb{Q}}
\newcommand{\R}{\mathbb{R}}

\newcommand{\SL}{\mathrm{SL}}

\newcommand{\sys}{\mathrm{sys}}

\newcommand{\tr}{\mathrm{tr}}

\newcommand{\Mg}{\mathcal{M}_g}

\def\tareesidedbox#1{\setbox0=\hbox{$#1$}\dimen0=\wd0 \advance\dimen0 by3pt\rlap{\hbox{\vrule height9pt width.4pt depth2pt \kern-.4pt\vrule height9.4pt width\dimen0 depth-9pt\kern-.4pt \vrule height9pt width.4pt depth2pt}} \relax \hbox to\dimen0{\hss$#1$\hss}}
\def\ho#1{\tareesidedbox{#1}}

\begin{document}

\title{Determining surfaces by short curves and applications}
%\author{Cayo D\'oria and Nara Paiva}

\author[Cayo D\'oria]{Cayo D\'oria}\thanks{C. D\'oria was partially supported by CNPq Grant 408834/2023-4}
\address{UFS, Departamento de Matem\'atica - Av. Marcelo D\'eda Chagas s/n, 49100-000. S\~ao Crist\'ov\~ao, Brazil}
\email{cayo@mat.ufs.br}

\author[Nara Paiva]{Nara Paiva}
\address{Instituto de Matem\'atica e Estat\'istica \\ Universidade Federal de Goi\'as \\ Goi\^ania - Brazil}
\email{nara\_reges\_paiva@discente.ufg.br}

\begin{abstract}
The goal of this work is to give new quantitative results about the distribution of semi-arithmetic hyperbolic surfaces in the moduli space of closed hyperbolic surfaces. We show that two coverings of genus \(g\) of a fixed arithmetic surface \(S\) are \(P(\frac{1}{g})\) apart from each other with respect to Teichmuller metric, where \(P\) is a polynomial depending only on \(S\) whose degree is universal. We also give a super-exponential upper bound for the number of semi-arithmetic hyperbolic surfaces with bounded genus, stretch and degree of the invariant trace field, generalizing for this class similar well known bounds for arithmetic hyperbolic surfaces. In order to get these results we establish, for any closed hyperbolic surface \(S\) with injectivity radius at least \(s\), a parametrization of the Teichmuller space by length functions whose values on \(S\) are bounded by a linear function (with constants depending only on \(s\)) on the logarithm of the genus of \(S.\) 
\end{abstract}

\maketitle

\section{Introduction}

In this article, a closed hyperbolic surface \(S\) means a closed surface equipped with a hyperbolic metric. By the Gauss-Bonnet theorem, \(S\) has genus \(g \geq 2\). Moreover, \(S\) is  geometrically equivalent to a quotient \(G \backslash \Hyp\), where \(\Hyp\) denotes the hyperbolic plane and \(G\) is a cocompact torsion-free discrete subgroup of the Lie group \(H=\SL(2,\R).\) Here we are mainly interested in the interplay between the hyperbolic geometry of the quotient and the arithmetic of the group \(G.\) 

For any integer \(g \geq 2\), let \(\Mg\) be the moduli space of closed hyperbolic surfaces of genus \(g\) and consider the respective Teichmuller space \(\mathcal{T}_g.\) 
%(Our first result is purely geometric)
It is known that for any \(g \geq 2\) there exists  \(6g-4\) marked simple curves on the surface of genus \(g\) such that any element in \(\mathcal{T}_g\) is uniquely determined by lengths of respective simple closed geodesics homotopic to these simples curves \cite{SS86}. It is a natural problem to find for any marked surface \(S \in \mathcal{T}_g\) a parametrization by lengths of closed geodesics whose values in \(S\) are the smallest possible. 

If we consider a marked surface \(S\) with injectivity radius \(s>0,\) it is possible to determine \(9g-9\) simple closed geodesics from a short pants decomposition of the surface which parameterize \(\mathcal{T}_g\) and such that these lengths of the geodesics on \(S\) are bounded by a linear function \(Ag+B\), where \(A\) is universal and \(B\) depends only on \(s\) (see \cite[Proposition 4.1.7]{Doria18}). 

Now we are able to improve the linear bound on \(g\) by a logarithmic bound. The main difference in this new construction is that we replace the classical pants decomposition by a decomposition in curves and arcs introduced by Parlier in \cite{Par18}. For a closed marked curve \(\beta\) in a marked hyperbolic surface \(S\), the length of the unique closed geodesic in \(S\) homotopic to \(\beta\) is denoted by \(L_{\beta}(S)\).  In Section \ref{parametrization} we will prove the following theorem:
\newpage

\begin{theorem}\label{main theorem intro}
Given $s>0 $, for any marked closed hyperbolic surface $S \in \mathcal{T}_g$ with injectivity radius at least \(s,\) there exist finitely many marked closed curves \(\gamma_1, \ldots, \gamma_n\) in \(S\) with \(n \leq 15g-15\) such that:
\begin{enumerate}
    \item \(L_{\gamma_i}(S) \leq 20\log(4g)+8~\mathrm{arcsinh}\left(\frac{1}{\sinh(\frac{s}{2})}\right)\) for all \(i=1,\ldots,n;\)
    \item The map \(\mathcal{T}_g \to \mathbb{R}^n\) whose coordinates are given by the corresponding length functions \(L_{\gamma_1},\ldots,L_{\gamma_n}\) is injective.
\end{enumerate}   
\end{theorem}

We remark that in some sense this result is optimal. Indeed, there exists a sequence of genera \(\{g_i\}\) with \(\lim\limits_{i \to \infty} g_i = \infty\) such that we can find a closed hyperbolic surface \(S_i\) of genus \(g_i\) with \(\ell(\gamma) \geq \frac{4}{3} \log(g_i)\) for all closed geodesic \(\gamma\) of \(S_i\) (see \cite[Theorem 1.10]{KSV07}). These examples imply that the upper bound of the theorem can be improved only in the constants.  Moreover, given a marked topological pants decomposition \(\mathcal{P}\) of a closed surface, for any set of marked closed curves which parameterize \(\mathcal{T}_g\), at least one of them, namely \(\gamma\), should intersect transversally some curve of \(\mathcal{P}\). Hence, if a marked hyperbolic surface \(S\) satisfies \(L_\eta(S) \leq \epsilon\) for all \(\eta \in \mathcal{P}\), then by the Collar Lemma (see \cite[Corollary 4.1.2]{buser2010geometry}) we have \(L_\gamma(S) \geq \mathrm{arcsinh}\left(\frac{1}{\sinh(\frac{\epsilon}{2})}\right),\) which can be very large for \(\epsilon\) sufficiently small. Thus we can not remove the dependence on $s$.

The theorem above will be our main tool in two quantitative results about arithmetic and semi-arithmetic hyperbolic surfaces. The precise definition of these objects will be given in Subsection \ref{semi-arithmetic}. 

Firstly,  we recall that if we fix a hyperbolic surface \(S\) of genus \(2\), for any \(g >2\) there exist finitely many points in \(\Mg\) which are coverings of \(S\). Hence there exists a minimal distance \(\delta_S(g)\), concerning Teichmuller metric, among all such coverings of genus \(g\). However, it follows from the proof of the Ehrenpreis conjecture given by Kahn-Markovic (\cite{KM15}) that \(\delta_S(g)\) goes to \(0\) when \(g\) grows. It does not follow directly from their proof in a quantitative way how this function goes to \(0\) for a general surface \(S\).

In Subsection \ref{minimal distance} we will combine the pivotal property of \emph{minimal spacing} (see Lemma \ref{Sarnak} or \cite{LS94}) of suitable arithmetic groups with Theorem \ref{main theorem intro} in order to show the following theorem, which says that for arithmetic surfaces the function \(\delta_S(g)\) is bounded below by a polynomial function on \(1/g\).  

\begin{theorem} \label{ minimal distance intro}
Let \(S\) be an arithmetic hyperbolic surface. There exists a constant \(A\) depending only on \(S\) such that the Teichmuller distance of any two  surfaces \(S',S'' \in \Mg\) which cover \(S\) is bounded below by \(\dfrac{A}{g^{240}}.\)    
\end{theorem}

We expect, by our proof of the theorem above and the sharpness of the Theorem \ref{main theorem intro}, that there is no a better bound than polynomial. Indeed, this would be true if the following question has positive answer which can be seen as a quantitative form of the Ehrenpreis conjecture.
\begin{question}
Do there exist, for any  hyperbolic surface \(S\), a polynomial \(Q\) depending only on \(S\) such that \(\delta_S(g) \leq Q(1/g)?\)
\end{question}

It follows from a result of Borel \cite{Borel81} that for any \(g \geq 2\) there exist only finitely many arithmetic hyperbolic surfaces in \(\Mg.\) Since arithmetic hyperbolic surfaces are rare, it is natural to generalize the notion of arithmeticity in such way that we can find a larger set of examples and that the number theory continues giving geometric results. Among all generalizations, the probably most studied has been the one given in \cite{SSW00} by Schaller and Wolfart, where it was defined \emph{semi-arithmetic} groups (see definition in Subsection \ref{semi-arithmetic}) by a slightly (but important) change in the definition of Fuchsian arithmetic groups. Very recently, various aspects of these groups have been explored in connections with multiple objects of mathematics (see for example \cite{CD22}, \cite{Gen12}, \cite{K15}, \cite{McMullen22}).

Once we have defined a semi-arithmetic hyperbolic surface, we have naturally a geometric invariant (the genus) and an arithmetic invariant (the degree of the invariant trace field). In \cite{BCDT23}, a new more intricate invariant is defined, namely \emph{stretch}, and the main result in [op.cit] implies that there are only finitely many examples of surfaces with these invariants bounded. We recall that any arithmetic hyperbolic surface is semi-arithmetic with stretch \(1\). 

It is possible to deduce from \cite{BGLS10} that for any \(d\) fixed, if \(\mathcal{AS}_{g,d}\) denotes the set of arithmetic hyperbolic surfaces of genus \(g\) with degree of the invariant trace field at most \(d,\) then 
\begin{equation}\label{asymptotic of restricted arithmetic locus}
    \lim_{g \to \infty} \frac{\log \#\mathcal{AS}_{g,d}}{g\log g}=2\
\end{equation}

The tools used in \cite{BGLS10} in order to obtain \eqref{asymptotic of restricted arithmetic locus} do not apply to non-arithmetic hyperbolic surfaces. Nevertheless, we are able to count semi-arithmetic hyperbolic surfaces using the arithmetic nature of the exponential of lengths of closed geodesics. Indeed, these numbers are always unit integers of some number field. Moreover, if we take the length of some closed geodesic, the degree of such field and the stretch give bound to the coefficients of the minimal polynomial of its exponential over \(\Q\). We will apply our first geometric theorem in Subsection \ref{counting subsection} for embedding a set of semi-arithmetic hyperbolic surfaces in a list of finite integers polynomials to obtain the following general estimate.

\begin{theorem} \label{counting intro}
For any $g \ge 2$, $d \ge 2$ and \(L \geq 1\), let \(\mathcal{SA}[g,d,L]\) be the set of semi-arithmetic hyperbolic surfaces of genus \(g\), with invariant trace field of degree at most \(d\) and stretch at most \(L\). Then, there exists a constant \(C>0\) depending only on \(d\) and \(L\) such that 
 \[ 2 \leq \limsup\limits_{g \to \infty}\frac{\log \# \mathcal{SA}[g,d,L]}{g \log(g)} \leq C.\]
 Furthermore, the constant \(C\) satisfies \(C \leq ULd^2 \) for some universal constant \(U>0.\)
\end{theorem}

\section{Preliminaries}

\subsection{The geometry of the moduli spaces}

In this subsection, we review the classical de\-fi\-ni\-tions of the Teichmuller and Moduli spaces of closed hyperbolic surfaces. Our references are \cite{buser2010geometry} and \cite{FM12}.

A \textit{closed hyperbolic surface} of genus $g\geq 2$ is a closed surface $S$ with $g$ handles equipped with a Riemannian metric of constant curvature \(-1\). For a closed hyperbolic surface \(S,\) its injectivity radius \(r>0\) it is equal to one half of its systole \(\sys(S)\), i.e. \(r=\frac{1}{2}\sys(S)=\frac{1}{2}\min\{\ell(\gamma)\}\) where \(\gamma\) runs through all closed geodesics in $S$. 

If we fix a closed surface \(R\) of genus \(g\), a \emph{marked} hyperbolic surface is a pair $(S,\phi)$ where \(\phi:R \to S\) is a homeomorphism. We say that two marked surfaces \((S,\phi)\) and \((S',\phi')\) are equivalent if there exists an isometry \(m:S \to S'\) such that \(\phi'\) and \(m \circ \phi\) are isotopic. The set of equivalence classes is known as \emph{Teichmuller space} and denoted by \(\mathcal{T}_g\).

It follows from two theorems of Teichmuller that for any pair \((S,\phi),(S',\phi')\) of marked hyperbolic surfaces, there exists a unique quasiconformal homeomorphism $h: S \to S'$ (see \cite[\S 11.1.2]{FM12} for quasiconformal maps definition) of minimum dilation $K_h$ among all homeomorphisms $f:S \to S'$ such that \(\phi'\) is homotopic to \(f \circ \phi\) (see \cite[Theorems 11.8 and 11.9]{FM12}). We define the \textit{Teichmuller distance} between the pair of marked surfaces \((S,\phi)\) and \((S',\phi')\) by \(\frac{1}{2}\log(K_h)\). 

The \textit{moduli space} $\mathcal{M}_g$ of hyperbolic surfaces of genus \(g\) is the set of equivalence classes of hyperbolic surfaces with respect to isometries. Let \(S,S'\) be representative points in \(\Mg\), we define the Teichmuller distance of the classes of \(S\) and \(S'\), denoted by \(d_T(S,S'),\) as the minimal of all distances of marked surfaces which represent the respective points (see \cite[Section 12.3]{FM12}). We note that we can suppose, replacing \(S'\) by an isometric copy of it if necessary, that the mapping \(h\) preserves marking, i.e. if \(d_T(S,S')\) is the Teichmuller distance between \((S,\phi)\) and \((S',\phi'),\) and \(h:S \to S'\) realizes the distance, we can suppose that  \((\phi')^{-1} \circ h \circ \phi: R \to R\) is isotopic to identity.

For any nontrivial free homotopy class \([\beta]\) of closed curve in \(R,\) and marked hyperbolic surface \((S,\phi)\), the length of the unique closed geodesic homotopic to \(\phi(\beta)\) with respect of the hyperbolic metric on \(S\) is denoted by \(L_{[\beta]}(S)\). For abuse notation, for any closed curve \(\beta \subset R\), and marked surface \(S\), \(L_{\beta}(S)\) means \(L_{[\beta]}(S).\) We will call \(L_\beta\) as a \emph{length function} in \(\mathcal{T}_g.\)  The length functions satisfy the following theorem (see \cite[Theorem 6.4.3]{buser2010geometry}).

\begin{theorem}\label{Teich distance}
Let $S, S'$ be closed hyperbolic surfaces of genus $g\geq 2$. If a homeomorphism $f: S \to
 S'$ is $K$-quasiconformal, then for each closed geodesic $\beta$ on $S$, the unique closed geodesic  $f[\beta]$ in the homotopy class of $f\circ\beta$ satisfies 
$$
\frac{1}{K} L_\beta(S)\leq L_{f[\beta]}(S')\leq K L_\beta(S).
$$
\end{theorem}

For future reference, we recall that by the classical Poincaré-Koebe uniformization theorem, every hyperbolic surface of genus at least \(2\) can be given as a quotient $G \backslash \mathbb{H}$, where $G$ is a torsion-free Fuchsian group. 

\subsection{Semi-Arithmetic hyperbolic surfaces} \label{semi-arithmetic}
In this section, we will give the definitions of arithmetic and semi-arithmetic surfaces. The main references for what follows are \cite{BCDT23}, \cite{Borel81}, \cite{MR03} and \cite{SSW00}. 

We start by defining what a semi-arithmetic group is. Let $G < \SL(2,\R)$ be a Fuchsian group, we say that $G$ is \textit{cofinite} if the quotient space $G \backslash \mathbb{H}$ has finite area. We will denote by $G^{(2)}$ the group generated by the squares of the elements of $G$, i.e. $G^{(2)}:=\langle T^2 \, | \, T \in G\rangle$. The \textit{trace set} of $G$ is the set $\mbox{Tr}(G):=\{|\tr(T)|\, |\, T \in G\}$. Finally, we define the \textit{invariant trace field} of $G$ to be $F=\mathbb{Q}(\mbox{Tr}(G^{(2)})) \subset \R$.  

A cofinite Fuchsian group $G$ is said to be \textit{semi-arithmetic} if the invariant trace field is a totally real number field and $\mbox{Tr}(G^{(2)})$ is contained in the ring of integers $\mathcal{O}_F$ of $F$. We note that any conjugated group and any finite index subgroup of a semi-arithmetic group is also semi-arithmetic with the same invariant trace field.

For $G$ to be called an \textit{arithmetic} group, it must be semi-arithmetic and satisfy the following condition: if $\phi$ is \(\Q\)-monomorphism (or simply an embedding) from $F$ to  $\mathbb{R}$, different from the identity, then $\phi(|\tr(T)|)\in [-2,2]$ for all $T \in G^{(2)}$.

Let \(G\) be a semi-arithmetic group with invariant trace field \(F.\) Consider the \(F-\) subalgebra of the algebra \(\mathrm{M}(2,\R)\) of \(2 \times 2\) real matrices, generated by \(G^{(2)}\). It turns out that this algebra is a quaternion algebra over \(F,\) namely the \emph{invariant quaternion algebra} of \(G\) and denoted by \(A(G).\) Moreover, \(G^{(2)}\) is contained in the group of elements of reduced norm one in an order \(\mathcal{O} \subset A(G)\), i.e. \(G^{(2)} \subset \mathcal{O}^1\). We also note that \(A(G)\) does not depend, up to \(F\) isomorphism, of any conjugated or finite index subgroup of \(G.\)

It is well established (see \cite{Borel81} for example) that there exists an \(F-\)algebra morphism \(\rho:A(G) \to \mathrm{M}(2,\R)^r \) for some \(r \leq [F:\Q]\) such that \(\rho(\mathcal{O}^1)\) is a cofinite discrete subgroup of isometries of the product \(\Hyp^r\). Since \(r\) depends only on \(A(G)\), we call it the \emph{arithmetic dimension} of \(G.\) This name comes from the fact that for the arithmetic Fuchsian group, we have \(r=1.\)

We define the \emph{stretch} of \(G\) as the infimum of all \(L>0\) such that there exists a \(L-\)Lipschitz map \(f:\Hyp \to \Hyp^r\), with respect to maximum metric on \(\Hyp^r,\) such that \(f(T \cdot z)=\rho(T) \cdot f(z)\) for all \(z \in \Hyp\) and \(T \in G^{(2)}.\) It follows from the nature of \(\rho\) that \(L \geq 1\). Moreover, for arithmetic groups or triangle groups, we have \(L=1.\) Hence, the arithmetic dimension and stretch measure how far a semi-arithmetic is from being arithmetic or triangular.

We say that a surface $S \in \mathcal{M}_g$ is semi-arithmetic (respectively arithmetic) if $S=G \backslash \Hyp$ for some semi-arithmetic (respectively arithmetic) Fuchsian group $G$ and we denote the set of such surfaces by \(\mathcal{SA}_{g}\) (respectively  \(\mathcal{A}_g\)). In particular, it is well defined the invariant trace field, invariant quaternion algebra (up to isomorphism), arithmetic dimension, and stretch of $S$ in a natural way since these invariants do not depend on the conjugacy class of a semi-arithmetic Fuchsian group.  

We conclude this subsection by recalling a relation between the geometry of a hyperbolic surface with its uniformizing group. Let $G$ be a Fuchsian group, for each $T \in G$ hyperbolic, i.e. with \(|\tr(T)|>2,\) we have a unique geodesic \(L_T \subset \Hyp\) fixed by \(T\) where \(T\) acts as a translation of displacement \(\ell(T)\) and hence projects on a closed geodesic of same length on the quotient \(G \backslash \Hyp\). Conjugating \(T\) by a suitable element in \(\SL(2,\R)\) we can suppose that \(L_T\) is the imaginary axis and then $e^{\ell(T)}=\lambda(T)$, where $\lambda(T)$ is the absolute value of the biggest eigenvalue of \(T^2\).

\subsection{Number theoretic facts about semi-arithmetic Fuchsian groups}

Let $\alpha$ be an algebraic integer and consider its minimal polynomial $P \in \mathbb{Z}[X]$. We say that a complex number $\beta$ is \emph{Galois conjugate} to $\alpha$ if $P(\beta)=0$. The \emph{house} of \(\alpha\) is the number \(\ho{\alpha}=\max\{|\beta| : \beta \mbox{ is Galois conjugate to } \alpha\}\). 

We say that an algebraic integer \(\alpha\) is \emph{reciprocal} if it is Galois conjugate to its inverse. For a reciprocal number \(\alpha\), the set of roots of its minimal polynomial \(P\) is invariant by the map \(x \to x^{-1}\). Hence \(P\) has even degree and satisfies \(P(X)=X^n P(1/X)\).

Now we will give an extension of a proposition already proved in \cite[Proposition 4.2]{BCDT23}. We will explore the arithmetic nature of eigenvalues of hyperbolic elements in semi-arithmetic Fuchsian groups and we will give a systolic estimate for semi-arithmetic hyperbolic surfaces. 

\begin{proposition}\label{systole lower bound and house bound} 
  Let \(S\) be a semi-arithmetic closed hyperbolic surface of invariant trace field of degree at most \(d\) and stretch at most \(L\), then:
  \begin{enumerate}
      \item For any closed geodesic \(\beta \subset S,\) the number $e^{\ell(\beta)}$ is a reciprocal algebraic integer of degree at most $2d$ and house at most \(e^{\ell(\beta) L}\).
      \item There exists some constant \(s>0\) which depends only on \(d\) and \(L\) such that \(S\) has injectivity radius at least \(s\).
  \end{enumerate}
\end{proposition}
\begin{proof}
We can write \(S=G \backslash \Hyp\) for some semi-arithmetic Fuchsian group.  As we have mentioned before, 
\begin{equation} \label{length-trace-eigenvalue} \ell(\beta)=\log(\lambda(T(\beta))) \mbox { for any closed geodesic } \beta \subset S,\end{equation} where \(T(\beta)\) is a hyperbolic element of \(G\) whose fixed geodesic projects on \(\beta\) and \(\lambda(T)\) denotes the biggest eigenvalue of \(T^2,\) for any hyperbolic element \(T.\)

Since  \(e^{\pm \ell(\beta)}\) are the roots of the polynomial \(X^{2}-\mathrm{tr}(T(\beta)^2)X+1 \in \mathcal{O}_F[X],\) we conclude that \(e^{\ell(\beta)}\) is a real reciprocal algebraic integer of degree bounded above by \(2d\). The estimate of the house of \(e^{\ell(\beta)}\) follows from \cite[Proposition 3.3]{BCDT23} and this proves \((1).\) 

For \((2)\) it is sufficient to give a lower bound for all closed geodesics of \(X\). By Dimitrov's proof of   Schinzel-Zassenhaus conjecture \cite[Theorem 1]{Dim19}, we have \(\ho{\scriptstyle{\lambda  (T(\beta))}}>2^{\frac{1}{4d}}\). Thus, for any closed geodesic \(\beta \subset S\) it holds
\[\ell(\beta) = \frac{\log(\lambda  (T(\beta))^L)}{L} \geq \frac{\log(\ho{{\scriptstyle \lambda  (T(\beta))}})}{L} \geq \frac{\log(2)}{4dL}.\]

\end{proof}

We finish this subsection with a lemma which counts the number of reciprocal algebraic integers with an upper bound on degree and house. Let \(\mathcal{U}_m(X)\) be the set of real reciprocal algebraic integers %with \(2r\) real Galois conjugates,
with degree at most \(2m\) and house bounded by \(X\).

\begin{lemma}\label{counting_by_house}
 Given positive integer $m$ and a number $X>1$, the set \(\mathcal{U}_m(X)\) is finite with at most $2m(4mX)^{m^2}$ elements.
\end{lemma}
\begin{proof}
 
  We can assigns for each \(\lambda \in \mathcal{U}_m(X)\) its minimal polynomial \(P_\lambda\). Since the preimage of any such polynomial has at most \(2m\) elements, the cardinality of \(\mathcal{U}_m(X)\) is bounded by \(2m\) times the number of minimal polynomials of elements in \(\mathcal{U}_m(X)\).
  
 Let \(P\) be the minimal polynomial of an element of \(\mathcal{U}_m(X)\), if \(P\) has degree \(2d\), \(d \leq m\), then \(P\) has at most \(d\) coefficients different of \(1\) since \(P\) is reciprocal. More precisely, \(P=X^{2d}+a_{1}X^{2d-1}+\ldots+a_dX^d+a_{d-1}X^{d-1}+\ldots+a_1X+1\), and for each \(i=1,\ldots,d\), $a_i$ is \((-1)^i\) times the $i$-th symmetric function on $2d$ variables applied in the roots of \(P\). Hence, $$|a_i| \leq {2d \choose i}(\ho{\alpha})^i \leq (2dX)^i.$$
 
Since each \(a_i\) is an integer, the number of possibilities for \(P\) with degree \(2d\) is bounded by
 \[2^{d}(2dX)^{\frac{d(d+1)}{2}}\]
 Hence, there exist at most \[\sum_{d=1}^m 2^{d}(2dX)^{\frac{d(d+1)}{2}} \leq  m2^m (2mX)^{\frac{m(m+1)}{2}}  \leq (4mX)^{m^2} \] possibilities for a minimal polynomial of an element of \(\mathcal{U}_m(X)\) (the last inequality can be verified applying \(\log\) in both sides).
\end{proof}

\begin{remark}
 In \cite{CH17}, Calegary and Huang prove more general estimates which give a better upper bound than in the lemma above. But for the application of this result in what follows their results do not give a substantial improvement. 
\end{remark}

\section{Length coordinates with logarithmic bound} \label{parametrization}

Let $S\in \mathcal{M}_g$ be a closed hyperbolic surface of genus $g\geq 2$. We denote by \(d\) the hyperbolic distance on \(S\) and for any set \(\mathcal{U}\) of closed geodesics or arcs in \(S\), \(\ell(\mathcal{U})\) denotes the set of their respective lengths \(\{\ell(u):u \in \mathcal{U}\}.\) 
We can decompose $S$ into curves $\Gamma=\{\gamma_1,\ldots, \gamma_k\}$ and arcs $\mathcal{A}=\{a_1,\ldots,a_{m}\}$ so that $S\setminus\{\Gamma, \mathcal{A}\}$ is a collection of geodesic right angled hexagons.

For each curve $\gamma \in \Gamma$ we fix the choice of two arcs that have one endpoint in $\gamma$, one arc on each side of $\gamma$. Let's call such endpoints $p_\gamma^-$ and $p_\gamma^+$. Thus, we can define the \textit{twist parameter} of $\gamma$, denoted by $\tau(\gamma)$, as the oriented distance of $p_\gamma^-$ and $p_\gamma^+$, ie, $|\tau (\gamma)|=d(p_\gamma^-,p_\gamma^+).$

The surface $S$ is uniquely determined by the twist parameters of each $\gamma\in\Gamma$, and by the lengths of the decomposition arcs, that is, by the sets $\tau(\Gamma)=\{\tau(\gamma_1),\ldots,\tau(\gamma_k)\}$ and $\ell(\mathcal{A})$ (\cite[Lemma 3.4]{Par18}). However, the surface can be determined just by the twist parameter of each $\gamma\in\Gamma$ and the length of some closed curves in $S$. For this, let's define a new class of curves.

Given $a\in \mathcal{A}$ and an orientation on $S\setminus \Gamma$, this defines an orientation on the boundary curves of $S\setminus \Gamma$ that are at the ends of $a$, say $\gamma_1$ and $\gamma_2$. Take $\gamma_a$ to be the closed geodesic in the free homotopy class of $\gamma_1\ast a\ast \gamma_2\ast a^{-1}.$ This curve will be called the \emph{chain} associated with $a$. Thus, the surface $S$ is uniquely determined by $\tau(\Gamma)$, $\ell(\Gamma)$ and $\ell(\Gamma_\mathcal{A})$ where $\Gamma_\mathcal{A}=\{\gamma_a\}_{a\in\mathcal{A}}$ (\cite[Lemma 3.5]{Par18}). We call the sets of curves $\Gamma$ and $\Gamma_\mathcal{A}$ the \textit{curve and chain system}.

The next proposition is an adapted result from \cite[Theorem 1.3]{Par18}. We will give a proof for reader convenience. For this, we recall two useful general lemmata (see \cite[Lemma 2.6 and Lemma 3.8]{Par18}). 
\begin{lemma}\label{minimal loop bound}
  For any \(S \in \mathcal{M}_g\)  and any \(x \in S,\) there exists a geodesic loop \(\delta_x\) based on \(x\) such that
  \[\ell(\delta_x) < 2 \cdot \mathrm{arccosh} \left(\frac{1}{2\sin \frac{\pi}{12g-6}}\right).\]
\end{lemma}
\begin{lemma}\label{distance to collar}
Let \(S \in \mathcal{M}_g\) and let \(c \subset S\) a simple geodesic loop. Consider \(\gamma\) be the unique simple closed geodesic freely homotopic to \(c\) and let \[\mathcal{C}(\gamma)=\{x \in S \mid d(x,\gamma) \leq w(\gamma)\}\] be the collar of \(\gamma\), where \(w(\gamma)=\mathrm{arcsinh} \left( \dfrac{1}{\sinh(\frac{\ell(\gamma)}{2})}\right) \). Then
\[\sup_{p \in c} \{d(p,\mathcal{C}(\gamma)\} < \log\left(\sinh\left( \frac{\ell(c)}{2}\right)\right).\]
\end{lemma}

We recall that for any \(s>0\), the \(s-\)\emph{thick part} of \(\Mg\) is the set \(\Mg^{\geq s}\) formed by hyperbolic surfaces with injectivity radius at least \(s\). 

\begin{proposition}\label{theo1}
Given \(\varepsilon>0\), any surface $S\in \mathcal{M}_{g}^{\geq \varepsilon}$ admits a curve and chain system $\Gamma$, $\Gamma_\mathcal{A}$ satisfying    
\begin{align*}
    \ell(\gamma) & < 2\log(4g)  \mbox{ for all } \gamma\in\Gamma \\
    \ell(\gamma_a) & < 8\log(4g)+4~ \mathrm{arcsinh} \left(\dfrac{1}{\sinh(\frac{\varepsilon}{2})}\right)  \mbox{ for all } a \in \mathcal{A} \\
    \ell(a) & < 2\log(4g)+2~ \mathrm{arcsinh} \left(\dfrac{1}{\sinh(\frac{\varepsilon}{2})}\right)  \mbox{ for all } a \in \mathcal{A}.
\end{align*}
Furthermore, this curve and chain system can be chosen such that for any \(\gamma \in \Gamma\) and side of \(\gamma\), there exists at least one arc \(a \in \mathcal{A}\) joining \(\gamma\) to \(\delta \in \Gamma\) with \(\delta \neq \gamma.\)
\end{proposition}
\begin{proof}
The proof follows the same circle of ideas of \cite[Theorem 1.3]{Par18}. For every \(x \in S\) we let $\delta_x$ be the shortest nontrivial loop based on \(x\). As the smallest loop, $\delta_x$ satisfies the Lemma \ref{minimal loop bound}
$$
\ell(\delta_x)< 2 \cdot \mathrm{arccosh} \left(\frac{1}{2\sin \frac{\pi}{12g-6}}\right) <2\log(4g).
$$

Now we consider the collection \(\tilde{\Gamma}\) of all simple closed geodesics in \(S\) which are free homotopic to \(\delta_x\) for some \(x \in S.\) Let \(\Gamma\) be a maximal disjoint set of \(\tilde{\Gamma}\).

To construct the set of arcs for the desired decomposition, we start by considering a decomposition in Voronoi cells around the boundaries of $S \setminus \Gamma$. That is, the surface $S$ is partitioned into sets
$$
V_\gamma=\{x\in S \mid d(x,\gamma)\leq d(x,\delta),\, \delta \; \mbox{boundary of} \; S\setminus \Gamma \}
$$
for each \(\gamma\) on the boundary of $S\setminus \Gamma$. Then, by Lemma \ref{distance to collar} and by straightforward manipulations, for any \(x\in V_\gamma\) we have
\begin{align*}
   d(x,\gamma) & < \log\left(\sinh\left( \frac{\ell(\delta_x)}{2}\right)\right)+\mathrm{arcsinh} \left( \dfrac{1}{\sinh(\frac{\ell(\gamma)}{2})}\right) \\
                   & < \log(4g) + \mathrm{arcsinh} \left(\dfrac{1}{\sinh(\frac{\varepsilon}{2})}\right)
\end{align*}

The points that belong to more than one cell are called the \textit{cut locus} of the decomposition. Such points form the edges when they belong to exactly two cells and the vertices when they belong to three or more cells. If the point belongs to exactly three cells, we call it \textit{trivalent vertex}. If the vertex is trivalent, for each edge of the cut locus connected to it, a geodesic arc is constructed between the boundaries of $S\setminus \Gamma$ associated with the vertex, crossing this edge. We include such arcs in the set $\mathcal{A}$ that we are constructing. For more details, see \cite{Par18}.

\begin{minipage}{0.50\textwidth}
    \begin{figure}[H]
    \centering
    \includegraphics[width=5.8cm]{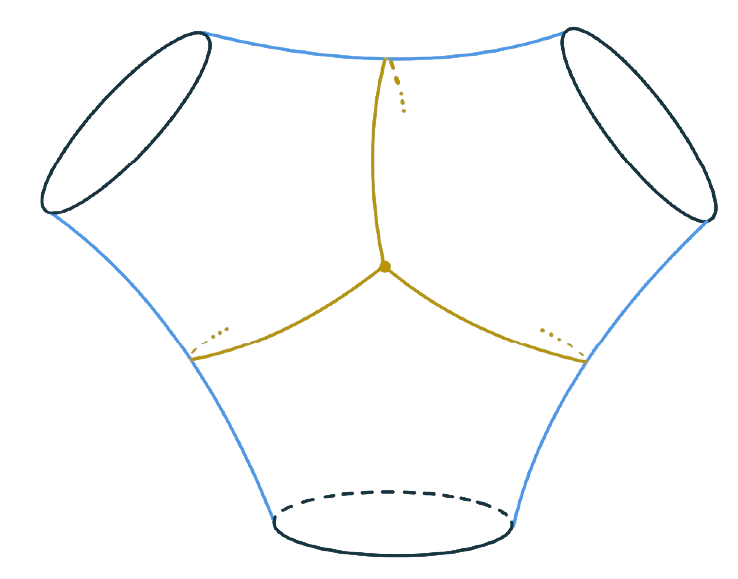}
    \caption*{Trivalent vertex}
    \label{fig1}
    \end{figure}
\end{minipage}
%\pause
\begin{minipage}{0.50\textwidth}
    \begin{figure}[H]
    \centering
    \includegraphics[width=5cm]{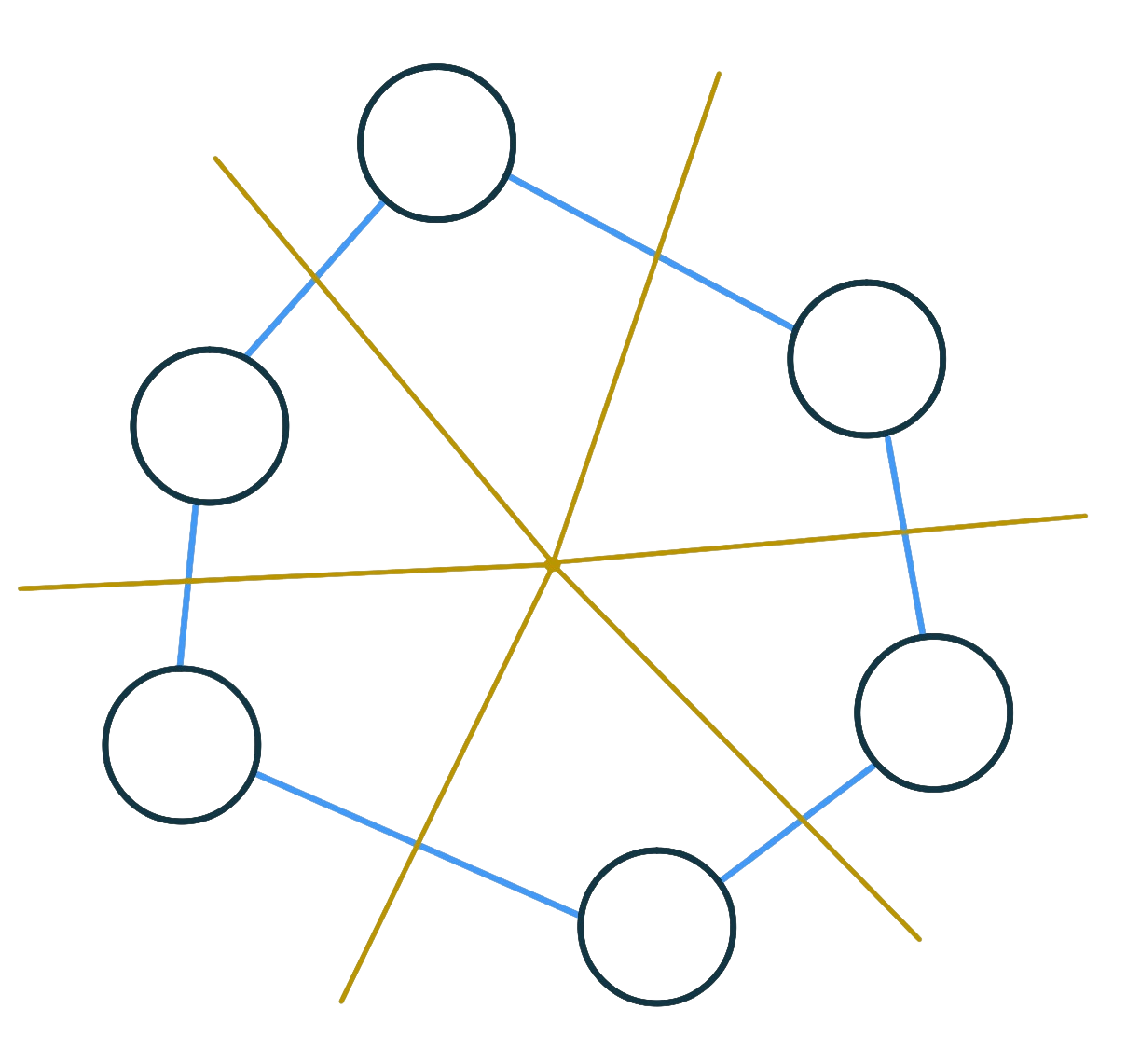}
    \caption*{Vertex with high degree}
    \label{fig2}
    \end{figure}
\end{minipage}
\vspace{0.2cm}

If the vertex $v$ has degree $q \geq 4$, we construct $q$ arcs as in the previous case. 

Given any of these $a\in \mathcal{A}$, let $\gamma'$ and $\gamma''$ be curves on the boundary of $S\setminus\Gamma$ joined by $a$. If $y$ is the point of intersection of the arc $a$ with the edge that separates the cells associated with $\gamma'$ and $\gamma''$, then $\ell(a)=2 d(y,\gamma')<2\log(4g)+2\, \mathrm{arcsinh}\left(\frac{1}{\sinh(\frac{\varepsilon}{2})}\right)$.
Now consider the $q$ simple arcs $\{c_i\}_{i=1}^q$ that connects the vertex $v$ to all associated boundaries of $S\setminus \Gamma$. We orient $c_1$ towards $v$ and all the other arcs orient from $v$. Now we concatenate $c_1$ with $c_i$, for $i=3,\ldots, q-1$, and take the arcs in its homotopy classes. So we add $q-3$ new arcs for set $\mathcal{A}$.

\begin{minipage}{0.34\textwidth}
    \begin{figure}[H]
    \centering
    \includegraphics[width=5cm]{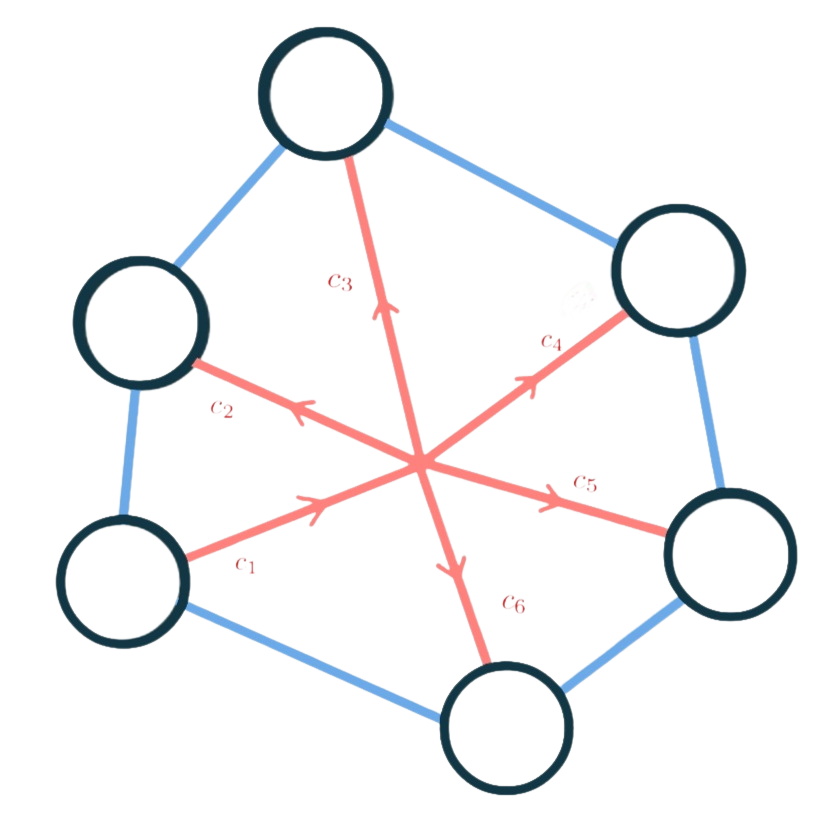}
    \caption*{\tiny Step 1: Arcs joining a vertex with boundary}
    \label{fig1}
    \end{figure}
    \vspace{0.1cm}
\end{minipage}
\begin{minipage}{0.34\textwidth}
    \begin{figure}[H]
    \centering
    \includegraphics[width=4.7cm]{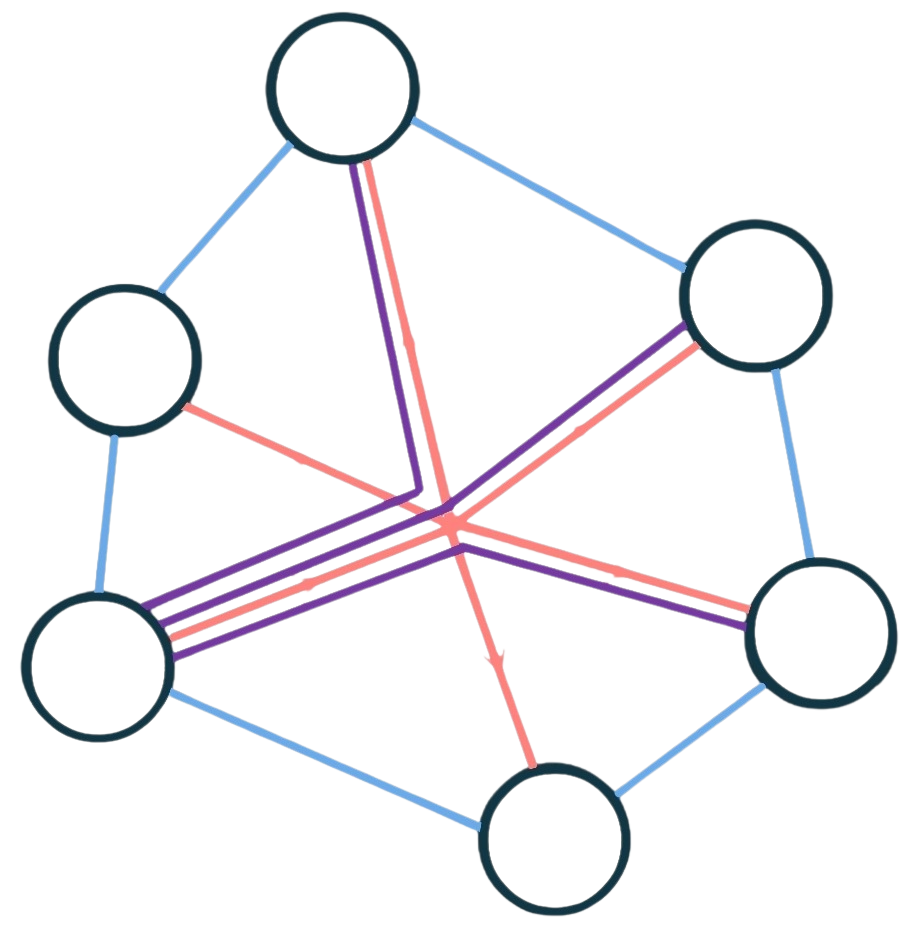}
    \caption*{\tiny Step 2: Concatenation of arcs \(c_i's\)}
    \label{fig2}
    \end{figure}
    \vspace{0.1cm}
\end{minipage}
\begin{minipage}{0.34\textwidth}
    \begin{figure}[H]
    \centering
    \includegraphics[width=5.2cm]{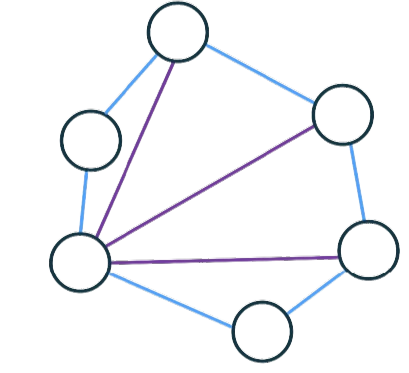}
    \caption*{\tiny Step 3: New arcs}
    \label{fig2}
    \end{figure}
    \vspace{0.15cm}
\end{minipage}

Hence, given an  arc \(a\) in the last construction with ends in $\gamma'$ and $\gamma''$, we have \(\ell(a) \leq 2~d(v,\gamma')\). Since \(\ell(\gamma_a) \leq 2\ell(a)+\ell(\gamma')+\ell(\gamma'')\), the inequalities above together prove the upper bound for all lengths \(\ell(\gamma),\, \ell(a)\) and \(\ell(\gamma_a).\)

In any case, the construction of the arcs is made so that all boundaries, in each component of \(S \setminus \Gamma\), are connected by at least two arcs. So after pasting, every geodesic in $\Gamma$ will have, for any side, at least one arc that connects it to a different geodesic in $\Gamma$. Therefore, we can choose arcs with this property in defining the twist parameter. 
\end{proof}

The proof of the Theorem \ref{main theorem intro} is a direct consequence of the following theorem.

\begin{theorem}\label{theo2}
Given $\varepsilon> 0$, let $S\in\mathcal{M}_{g}^{\geq \varepsilon}$. There exists a curve and chain system $\Gamma$, $\Gamma_\mathcal{A}$ for the surface $S$ and an associated set of closed geodesics $\Lambda=\{\delta^\alpha,\eta^\alpha\}_{\alpha\in\tau(\Gamma)}$ such that $S$ is uniquely determined by $\ell(\Gamma)$, $\ell(\Gamma_\mathcal{A})$ and $\ell(\Lambda)$. The length of the curves in $\Gamma$, $\Gamma_\mathcal{A}$ and $\Lambda$ and also the length of the arcs in $\mathcal{A}$ are bounded above by
$$
 20\log(4g)+8~ \mathrm{arcsinh}\left(\frac{1}{\sinh(\frac{\varepsilon}{2})}\right).
$$
Furthermore, the amount of curves that determine the surface $S$ is at most $15g-15$.
\end{theorem}
\begin{proof}
Consider the construction of the curve and chain system $\Gamma, \Gamma_\mathcal{A}$ from the previous proposition. Given $\gamma \in \Gamma$, by Proposition \ref{theo1} we can choose arcs $a_1$ and $a_2$ on $\mathcal{A}$ connected to $\gamma$, on opposite sides, such that \(a_i\) connect \(\gamma\) to $\gamma_i \in \Gamma$ with $\gamma_i \neq \gamma$ for \(i=1,2\). The surfaces with boundaries $\gamma$, $\gamma_i$ and $\gamma_{a_i}$ are Y-pieces denoted by $Y_i$, with $i=1,2$, embedded in $S$. 

Looking at surface $S\setminus \{\gamma\}$, let's reparameterize its boundaries $\bar{\gamma}_1$ and $\bar{\gamma}_2$ so that
$$
\bar{\gamma}_1(0)=p_\gamma^- \quad \text{and} \quad \bar{\gamma}_2(0)=p_\gamma^+.
$$
Let $\alpha=\tau(\gamma)$, since
$$
\tau(\gamma)=\pm d(p_\gamma^-,p_\gamma^+)=\pm d(\bar{\gamma}_1(0),\bar{\gamma}_2(0))
$$
we can reparameterize $\gamma$ in $S$, by arc length, and rename it $\gamma^\alpha$, where
\begin{equation}\label{eq1}
\gamma^\alpha(t):=\bar{\gamma}_1(t)=\bar{\gamma}_2(\alpha-t), \quad t\in \R / [t \to t +\ell(\gamma)]
\end{equation}

Hence, we have an $X$-piece immersed in $S$ defined by
$$
X^\alpha :=Y_1+Y_2\;\; \text{mod}(\ref{eq1}).
$$

Let $\pi_\alpha: X^\alpha\rightarrow S$ be the natural projection. We define $\delta^\alpha$ and $\eta^\alpha$ the closed geodesics in the homotopy class of $\pi_\alpha(\tilde{\delta}^\alpha)$ and $\pi_\alpha(\tilde{\eta}^\alpha)$ respectively, where $\tilde{\delta}^\alpha$ and $\tilde{\eta}^\alpha$ are curves in $X^\alpha$ that we define below.

Since $X^\alpha$ and $S$ have the same twist parameter on \(\gamma\) and since the projection $\pi_\alpha$ is length-preserving, by [Proposition 3.3.14, \cite{buser2010geometry}] there is an analytic function $F$ satisfying
$$
\alpha= F(\ell(\gamma^\alpha), \ell(\delta^\alpha), \ell(\eta^\alpha), \ell(\gamma_1), \ell(\gamma_2), \ell(\gamma_{a_1}), \ell(\gamma_{a_2})).
$$ 
Let $b$ to be the arc parameterized in $\gamma^\alpha$ that realizes the distance of $p_\gamma^-$ and $p_\gamma^+$. And consider $d^\alpha$ be the arc, perpendicular to $\gamma_1$ and $\gamma_2$, in the homotopy class of $a_2\ast b\ast a_1^{-1}$. 

    \begin{figure}[H]
    \centering
    \includegraphics[width=5cm]{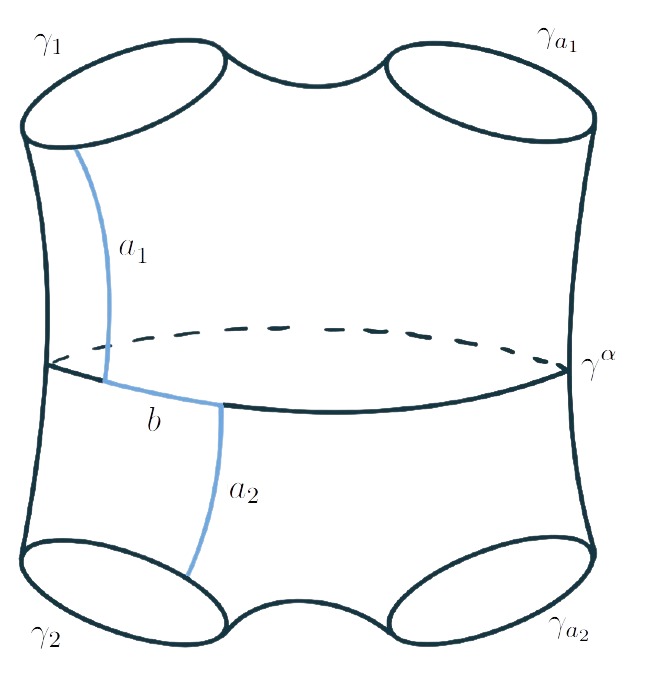}
    \end{figure}
The curve $\tilde\delta^\alpha$ is defined as the closed geodesic in the homotopy class of $d^\alpha\ast \gamma_1\ast (d^\alpha)^{-1}\ast \gamma_2$. And the curve $\tilde\eta^\alpha$ is homotopic to the image of the twist of $\tilde\delta^\alpha$ in $\tilde\gamma^\alpha$ for a complete turn (for more details see \cite[Subsection 3.3]{buser2010geometry}).

Therefore,
\begin{eqnarray*}
    \ell(\delta^\alpha)&< & \ell(\gamma_1)+\ell(\gamma_2)+2\ell(\gamma^\alpha)+2\ell(a_1)+2\ell(a_2)\\
    &< & 16\log(4g)+8~ \mathrm{arcsinh}\left(\frac{1}{\sinh(\frac{\varepsilon}{2})}\right)
\end{eqnarray*}

Now, the curve $\eta^\alpha$ has the same length as curve $\delta^{\alpha+\ell(\gamma^\alpha)}$, so
$$
\ell(\eta^\alpha)<20\log(4g)+8~ \mathrm{arcsinh}\left(\frac{1}{\sinh(\frac{\varepsilon}{2})}\right).
$$

We conclude that the surface $S$ is determined by $\ell(\Gamma)$, $\ell(\Lambda)$ and $\ell(\Gamma_\mathcal{A})$.
Note that in the decomposition we obtain $4g-4$ hexagons since the surface area $S$ is $4\pi(g-1)$ and the area of each hexagon is $\pi$. In each hexagon, we count $3$ arcs. Then we have $12g-12$ arcs, which becomes $6g-6$ arcs in $\mathcal{A}$ after gluing. Since there are most \(3g-3\) simple closed geodesics disjoint in \(S\), \(\Gamma\) has \(k \leq 3g-3\) elements. As each $\delta^\alpha$ and each $\eta^\alpha$ are associated with a $\gamma\in\Gamma$, $\#\Lambda = 2k$. Thus we have that $S$ is determined by \[\#\Gamma +\#\Lambda +\# \Gamma_\mathcal{A}\leq 3k+6g+6 \leq 15g-15\] curves. 
\end{proof}

We finish this section with a lemma which will be useful in Subsection \ref{counting subsection}. We say that a family \(\Omega\) of closed curves on a surface \(R\) gives a parametrization by a marked (topological) curve and chain system if \(\Omega\) is obtained by a construction as in the proof of the theorem above. By definition \(\Omega\) is constructed by a unique curve and chain system, which we will call the \emph{underlying curve and chain system} of \(\Omega.\)

Let \(\mathcal{F}_g\) be the set of all possible \(\Omega's\) for a surface of genus \(g\) up to the action of the group of homeomorphisms of \(R\). Since this action on \(\Mg\) is trivial, then does not matter if we replace a curve and chain system on a surface by its image by a homeomorphism. But, topologically, the set of all systems can be decomposed into finitely many orbits. Now we are able to give a quantitative upper bound for \(\# \mathcal{F}_g.\)

\begin{lemma}\label{topological types}
   The set \(\mathcal{F}_g\) is finite and 
   \[\# \mathcal{F}_g \leq B^g g^{6g} \]
for some universal constant \(B>0.\) 
\end{lemma}
\begin{proof}
Consider the map which associates for any \(\Omega \in \mathcal{F}_g\) its underlying curve and chain system. If this system has \(k\) closed curves, we get at most \(k \cdot (\#\mathcal{A})^2\) possibilities for choosing chains in order to recover \(\Omega\). Hence, if we denote by \(N_{cc}(g)\) the number of orbits of curves and chain systems on a surface of genus \(g\)  by the action of orientation-preserving homeomorphism group, then \(\# \mathcal{F}_g\) is bounded above by \((3g-3)(4g-4)^2 \cdot N_{cc}(g).\) 
 
By \cite[Lemma 3.2]{Par18} we have the following estimate
\[N_{cc}(g) \leq \frac{1}{e^6}\left(\frac{12^6}{e^5} \right)^{g-1}(g-1)^{6g-6}. \]
Therefore,
\[\# \mathcal{F}_g \leq \frac{48(g-1)^2}{e^6}\left(\frac{12^6}{e^5} \right)^{g-1}(g-1)^{6g-6} \leq B^g g^{6g}\]
for some universal constant \(B>0.\)

\end{proof}

\section{Applications}
\subsection{Minimal distance between coverings of arithmetic hyperbolic surfaces}\label{minimal distance}

Given a hyperbolic surface \(S\), for any \(n>1\) the number of \(n\)-sheeted coverings of \(S\) is equal to the number of conjugacy classes of subgroups of index \(n\) of the Fuchsian group \(\Gamma\) which uniformizes \(S\). On the other hand, for a finitely generated group, there exist finitely many subgroups of a given index. By the Gauss-Bonnet theorem, the set of points in \(\Mg\) which are coverings of \(S\) is equal to the number of \(m\)-sheeted coverings of \(S\) for some explicit \(m\) depending only \(g\) and the genus of \(S.\) Thus, it is well defined the function \(\delta_S\) given by 
\begin{align*}
  \delta_S(g) =    \begin{cases}
   0 , \mbox{ if there not exists coverings of } S \mbox{ of genus } g \\
    \min\{d_T(\tilde{S},\overline{S})~:~ \tilde{S},\overline{S} \in \Mg \mbox{ are different coverings of S} \}, \mbox{ otherwise,}
    \end{cases}
\end{align*}
where \(d_T\) denotes the Teichmuller metric in \(\Mg\). 

It follows from the work of Kahn and Markovic in \cite{KM15}, where they proved the Ehrenpreis conjecture, that for any \(\varepsilon>0\) there exist different coverings of a fixed surface \(S\) with Teichmuller distance at most \(\varepsilon.\) Hence, \(\lim\limits_{g \to \infty} \delta_S(g) = 0 \) for any \(S.\)

It would be interesting if there would exist a uniform (qualitative) upper bound of all such functions. For example, is there a polynomial \(Q\) depending only on \(S\) such that \(\delta_S(g) \leq Q(1/g)?\)
If the answer is positive then this is the better which we can obtain due to our next theorem. Before it, we will recall the following Lemma - originally proved in \cite[Lemma 2.1]{LS94} - whose proof is simple and captures a strong particularity of arithmetic surfaces.

\begin{lemma} \label{Sarnak}
If $G$ is an arithmetic Fuchsian group with invariant trace field $F,$ then there exists a constant $c>0$ which depends only on $F$, such that for any $T, T' \in G$ with $|\mathrm{tr}(T)| \neq |\mathrm{tr}(T')|$ we have $$ |\mathrm{tr}(T)^2-\mathrm{tr}(T')^2| > c. $$
\end{lemma}
\begin{proof}
  If $F$ is the invariant trace field of $G$ and $\phi$ is any nontrivial embedding of this field into $\mathbb{C},$ then for any trace $t \in \{ \mathrm{tr}(T^2) \, | \, T \in G \}$ we have $|\phi(t)| \leq 2$. Since \(\tr(A^2)=\tr(A)^2-2\) for any \(A \in \SL(2,\R)\), if $|\mathrm{tr}(T)| \neq |\mathrm{tr}(T')|,$ then \(\tr(T^2)=\tr(T)^2-2 \neq \tr(T')^2-2=\tr(T'^2).\) 
 
 Let \(\prod\limits_{i=1}^d \phi_i(x)\) be the \emph{norm} of any \(x \in F\), where \(\{\phi_i\}\) denotes the set of embedding of \(F\) in \(\R\) with \(\phi_1=\mathrm{id}\) and \(d=[F:\Q]\) is the degree of \(F\). For \(\xi \in \mathcal{O}_F,\) its norm is a nonzero integer. Thus,
 $$ 1 \leq |\mathrm{tr}(T^2)-\mathrm{tr}(T'^2)| \prod\limits_{i=2}^d |\phi_i(\mathrm{tr}(T^2))-\phi_i(\mathrm{tr}(T'^2)) | \leq   |\mathrm{tr}(T)^2-\mathrm{tr}(T')^2| 4^{d-1} .$$
\end{proof}

For any closed geodesic $\beta$ on a hyperbolic surface \(S=\Lambda \backslash \Hyp\), let $T_\beta(S) \in \Lambda$ be a hyperbolic element such that its axis projects on $\beta$. It follows from \eqref{length-trace-eigenvalue} that the length $\ell_{\beta}(S)$  and $\mathrm{tr}(T_\beta(S))$ are related by the equality:
\begin{equation} \label{traceA}
\ell_{\beta}(S)=2\cosh^{-1}\left(\frac{\tr(T_\beta(S))}{2}\right).
\end{equation}

Theorem \ref{ minimal distance intro} follows directly from the following more general theorem.
\begin{theorem}
 Let $S$ be an arithmetic closed surface of genus $g \geq 2$ with invariant trace field \(F\). Then there exist constant $A=A(d)$, which depends only on the degree \(d\) of \(F,\) such that if  $S' \in \mathcal{M}_g$ is another arithmetic closed surface with invariant trace field \(F\), it holds
 $$
 d_T(S,S') \geq A\, g^{-240}.
 $$
\end{theorem}

\begin{proof}
We fix $S=G \backslash \mathbb{H} \in \mathcal{M}_g$ an arithmetic closed geodesic with invariant trace field \(F\) of degree \(d\). Since any arithmetic surface has stretch \(1,\) by Proposition \ref{systole lower bound and house bound}, there exists \(s>0\) depending only on \(d\) such that \(S \in \mathcal{M}_{g}^{\geq s}.\) Now we can apply the Theorem \ref{theo2} in order to get a curve and chain system $\Gamma=\{ \gamma_1, \ldots, \gamma_k\}$, $\Gamma_\mathcal{A}=\{\gamma_a \mid a\in\mathcal{A}\}$ of $S$ such that $S$ is completely determined by $\Omega=\{\gamma_i,\delta_i, \eta_i,\gamma_{a} \mid 1\leq i\leq k,\; a \in \mathcal{A}\}$ and 
$$
\ell(\gamma_i), \ell(\delta_i), \ell(\eta_i), \ell(\gamma_{a})  \leq  60\log g + c(d) \hspace{0.3cm} \mbox{for all} \hspace{0.3cm} 1 \leq i \leq k  \mbox{ and } a \in \mathcal{A}.
$$ 
for some explicit constant \(c(d)>0.\)

Let $S' \neq S \in \mathcal{M}_g$ be an arithmetic surface with invariant trace field \(F\). 
If \(d_T(S,S') \geq \frac{\log(2)}{2}\) then the theorem is proved with \(A=1.\) Thus we can suppose  $d_T(S,S')<\frac{\log(2)}{2}$. 
Let  $f:S \to S'$ be a $\mathrm{e}^{2\varepsilon}$-quasiconformal map for some \(\varepsilon<\frac{\log(2)}{2}\) which we can suppose isotopic to identity.
Since $S\neq S'$, for some $\omega \in \Omega $ we have $\ell_{\omega}(S) \neq \ell_{\omega}(S')$. Let $G'$ be the uniformizing group of \(S'\). If we denote by \(t=\tr(T_\omega(S))\) and \(t'=\tr(T_\omega(S'))\), then \(|t| \neq |t'|\) and by Lemma \ref{Sarnak} it holds
\begin{equation} \label{traceB} 
|t^2-(t')^2| > c
\end{equation}
for some constant \(c>0\) which also depends only on \(d.\)

By Theorem \ref{Teich distance}, we obtain the following relation
$$
\mathrm{e}^{-2\varepsilon}\leq \frac{\cosh^{-1}\left(\frac{(t')^2-2}{2}\right)}{\cosh^{-1}\left(\frac{t^2-2}{2}\right)} \leq \mathrm{e}^{2\varepsilon}
$$
Hence,
$$
\left|\log\cosh^{-1}\left(\frac{(t')^2-2}{2}\right)-\log\cosh^{-1}\left(\frac{t^2-2}{2}\right)\right|< 2\,\varepsilon.
$$

Let $L(x)=\log\cosh^{-1}\left(\frac{x-2}{2} \right)$. By the mean value theorem, we know that there exists $\theta$ between $t^2$ and $(t')^2$ such that
\begin{equation}\label{nara1}
c\,L'(\theta)<|(t')^2-t^2|\,|L'(\theta)|=|L((t')^2)-L(t^2)|<2\,\varepsilon.
\end{equation}

Now we want to estimate \(L'(\theta)\) by below. Since \(\log(x) \leq \cosh^{-1}(x)\) for any \(x>1\), we have 
\begin{align*}
2\log \left(\frac{|t|}{2}\right) & \leq  2\cosh^{-1}\left(\frac{|t|}{2}\right) =  \ell_\omega(S)\leq 60\log g + c(d)=\log(C_d \cdot g^{60})
\end{align*}
where \(C_d=e^{c(d)}\).
Moreover,
\begin{align*}
 2\, \log\left(\frac{|t'|}{2}\right)& \leq  2\cosh^{-1}\left(\frac{|t'|}{2}\right) = \ell_\omega(S')\leq \mathrm{e}^{2\varepsilon} \cdot L_\omega(S) 
\end{align*}

We recall that we are assuming $\mathrm{e}^{2\varepsilon}<2$.  Thus, \begin{align*}t^2 \leq C_d \cdot 4 \cdot g^{60} \mbox{ and  } (t')^2 \leq C^2_d\cdot 4 \cdot g^{120}     
\end{align*} 
Hence, \(4< \theta \leq \max\{t^2,(t')^2\}  \leq 4\cdot C^2_d g^{120}\).
By a direct calculation, we get for any \(x>4\)
$$
L'(x)=\frac{1}{\sqrt{x(x-4)}  \, \cosh^{-1}\left(\frac{x-2}{2} \right)}>\frac{1}{x^2}.
$$
Therefore, \(L'(\theta)>\dfrac{1}{\theta^2} >\dfrac{1}{16~C_d^4~g^{240}}\). 

Returning to the equation \eqref{nara1}, we obtain $\varepsilon>A\, g^{-240},$ where $A=\frac{c}{16~C_d^4}$ depends only on \(d\).
\end{proof}
 We remark that the degree \(240\) is far from sharp. Indeed, it depends on non-optimal constants given by Theorem \ref{theo2}. Moreover, in the argument of the proof above we could improve the upper bound of \((t')^2\) if we assume that \(d(S,S') < \varepsilon\) for any \(\varepsilon>0\) given, which would increase the constant \(A\) at the end. 

\subsection{Counting semi-arithmetic hyperbolic surfaces}\label{counting subsection}

In the previous subsection we have used the fact that for any given totally real number field \(F\), there exist finitely many arithmetic hyperbolic surfaces with invariant trace field \(F\) of a fixed genus. This fact is no longer true for semi-arithmetic surfaces. Indeed, in \cite{CD22} it is shown that for any \(g \geq 2\) there exists a family of semi-arithmetic hyperbolic surfaces of genus \(g\) with the same invariant trace field and arbitrarily small injectivity radius. However, if we bound the arithmetic dimension and stretch we have a picture similar to the arithmetic case. 

For any integer \(r \geq 1\) and positive number \(L \geq 1\), we define the set 
\[\mathcal{SA}_{g}(r,L)=\{ S \in \mathcal{SA}_{g} \mid  \mathrm{a.dim}(S) \leq r \mbox{ and }  \delta(S) \leq L\},\]
where \(\mathrm{a.dim}(S)\) is the arithmetic dimension of \(S\) and \(\delta(S)\) is the stretch of \(S\). In \cite{BCDT23}, it is shown that \(\mathcal{SA}_{g}(r,L)\) is always finite. In particular, for any \(d \geq 1\), there exist finitely many semi-arithmetic closed hyperbolic surfaces of genus \(g\) with bounded stretch and invariant trace field of degree at most \(d\), since the arithmetic dimension is bounded above by the degree of the invariant trace field. 

We define now the set \(\mathcal{SA}[g,d,L]\) of semi-arithmetic closed hyperbolic surfaces of genus \(g\) with invariant trace field of degree at most \(d\) and stretch at most \(L\). We note that the set of arithmetic hyperbolic surfaces of genus \(g\) with invariant trace field of degree at most \(d\), denoted by \(\mathcal{AS}_{g,d}\), is always a subset of \(\mathcal{SA}[g,d,L]\) for any \(L \geq 1.\)

Now we fix the well known Bolza surface \cite{KKSV16}, which is an arithmetic hyperbolic surface of genus \(2\) with quadratic invariant trace field. Hence, the number of coverings with genus \(g\) of the Bolza surface gives a lower bound for \(\mathcal{SA}[g,d,L]\) for any \(L \geq 1\) and \(d \geq 2\). We note that it can be traced back from the proof of \cite[Corollary 1.4]{BGLS10} that
\begin{equation}\label{lower bounding for counting}
\#  \mathcal{SA}[g,d,L] \geq (bg)^{2g}   
\end{equation}
for some universal constant \(b>0.\) 

We can apply our main theorem in order to show that this type of grow is sharp. 

\begin{theorem} \label{thm:surfaces}
 For any $g \ge 2$, $d \ge 2$ and \(L \geq 1\), there exists a constant \(C>0\) depending only on \(d\) and \(L\) such that 
 \[ 2 \leq \limsup\limits_{g \to \infty} \frac{\log(\# \mathcal{SA}[g,d,L])}{g \log(g)} \leq C.\]
 Furthermore, the constant \(C\) satisfies \(C \leq ULd^2 \) for some universal constant \(U \geq \frac{1}{2}.\)
\end{theorem}

The idea of the proof of this theorem is to use arithmetic information of hyperbolic elements in a semi-arithmetic Fuchsian group by its invariants and apply the main theorem to reducing the problem into counting real reciprocal algebraic integers with bounded degree and house.

\begin{proof}[Proof of the Theorem \ref{thm:surfaces}]
The lower bound follows from \eqref{lower bounding for counting}. For the upper bound, let $g>1, d \geq 2$ and $L \geq 1$ be fixed parameters. 
For each \(S^* \in \mathcal{SA}[g,d,L]\), by Proposition \ref{systole lower bound and house bound} there exists a constant \(s=s(d,L)\) such that \(S^* \in \Mg^{ \geq s}\). By Theorem \ref{main theorem intro} there exists $n \leq 15(g-1)$ length functions $L_1,\ldots,L_n,$ such that \(\{L_i(\cdot)\}_i\) parametrizes \( \mathcal{T}_g\) and $L_i(S) \leq 20\log(\sigma g)$ for all $i=1,\ldots, k$, where \(\sigma\) is a constant which depends only on \(d\) and \(L\), where \(S\) is some representative of \(S^*\) in \(\mathcal{T}_g\). 

By Proposition \ref{systole lower bound and house bound}, each $\exp(L_i(S))$ is an algebraic integer of degree at most $2d$ and \[\ho{\exp{\scriptstyle (L_i(S))}} \leq \exp(L_i(S)\cdot L)\leq \exp(\log(\sigma g) \cdot 20L)=(\sigma g)^{20L}.\]
 
By Lemma \ref{counting_by_house}, there are at most $$2d\left[4d \cdot (\sigma g)^{20L}\right]^{d^2}$$  possibilities for each entry of $P(S)$.

If $P(S)=P(S')$ with $S \neq S'\) where \(S' \in \mathcal{T}_g\) projects onto some $S'^{*} \in \mathcal{SA}[g,d,L]$, then necessarily the corresponding lists come from different topological types of parametrization by marked curve and chain decomposition. By Lemma \ref{topological types}, the number of topological types of such parametrization in a closed surface of genus $g$ is at most $B^g g^{6g}$ for some universal constant $B>0$. Hence, 
$$\# \mathcal{SA}[g,d,L] \leq B^{g}g^{6g}\left[2d(4d)^{d^2}\cdot (\sigma g)^{20Ld^2}\right]^{15(g-1)} \leq B_1^g \cdot g^{B_2 g}$$
for constants \(B_1,B_2>0\), with \[B_1=B_1(d,L) \leq B[2d(4d)^{d^2}\sigma^{20Ld^2}]^{15} \mbox{ and } B_2=B_2(d,L) \leq 6+300Ld^2.\]
We put all above inequalities together and we obtain
\[\frac{\log(\# \mathcal{SA}[g,d,L])}{g\log(g)} \leq \frac{\log(B_1)}{\log(g)}+B_2 \leq C\]
where \(C=C(d,L)=\frac{\log(B_1(d,L))}{\log(2)}+B_2(d,L)\) depends only on \(d\) and \(L\) in such way that 
\(C(d,L) \leq U Ld^2\) for some universal constant \(U>0.\)
\end{proof}
 We recall that the notation \(f(x)=O_{a,b}(h(x))\), for positive real functions \(f,h\), means that there exist \(x_0>0\) and a constant \(M>0\) depending only on \(a,b\) such that \(f(x) \leq Mg(x)\) for all \(x>x_0.\) 
 
\begin{corollary}
For any \(r,L \geq 1\) fixed, it holds
\[\log(\# \mathcal{SA}_g(r,L)) = O_{r,L}(g^{1+\varepsilon}) \mbox{ for every } \varepsilon>0 \mbox{ when } g \to \infty\]
\end{corollary}
\begin{proof}
When we fix \(r,L\), the theorem \cite[Theorem 4.1]{BCDT23} tells us that for any \(S \in \mathcal{SA}_g(r,L) \), its invariant trace field has degree at most \(C' \log(g)\) for some constant \(C'\) which depends only on \(r\) and \(L.\) The corollary follows from this inequality, the estimate of the constant \(C\) in Theorem \ref{thm:surfaces} and the fact that for any natural \(n\) we have \(\log^n(x)=O(x^{\varepsilon})\) for every \(\varepsilon>0.\)
\end{proof}
\bibliographystyle{plain}

\bibliography{ref}
\end{document}